\def\timestamp{%
Time-stamp: <codiag-notconfl.tex: Saturday 11-10-2008 at 16:14:12 (cest)>}
\def\stripname Time-stamp: <#1 #2>{#2}
\edef\filedate{\expandafter\stripname\timestamp}

\documentclass[a4paper]{amsart}

\usepackage{amssymb,amsrefs}

\newcommand{\cl}{\operatorname{cl}}
\newcommand{\orpr}[2]{\langle#1,#2\rangle}
\newcommand{\preim}{^\gets}

\DeclareMathSymbol\restr2{AMSa}{"16}

\newtheorem{claim}{Claim}

\DeclareSymbolFont{cmmib}{OML}{cmm}{b}{it}
\DeclareMathSymbol\xv0{cmmib}{`x}

\theoremstyle{remark}
\newtheorem*{remark}{Remark}

\begin{document}

\title{A concrete co-existential map that is not confluent}

\author{Klaas Pieter Hart}
\address{Faculty of Electrical Engineering, Mathematics and Computer Science\\
         TU Delft\\
         Postbus 5031\\
         2600~GA {} Delft\\
         the Netherlands}
\email{k.p.hart@tudelft.nl}
\urladdr{http://fa.its.tudelft.nl/\~{}hart}

\date{\filedate}

\begin{abstract}
We give a concrete example of a co-existential map between continua 
that is not confluent.
\end{abstract}


\maketitle

\section*{Introduction}

In \cite{bankston-conf} Bankston gave an example of a co-existential map that 
is not confluent.
The construction is rather involved and does not produce a concrete example
of such a map.
A lot of effort was needed to get the main ingredient, to wit a co-diagonal 
map that is not monotone.

The purpose of this note is to show that one can write down a concrete map
between two rather simple continua that is co-existential and not confluent.
It will be clear from the construction that the range space admits co-diagonal
maps that are nor confluent and, a fortiori, not monotone.

\section{Preliminaries}

In the interest of brevity we try to keep the notation down to 
the bare minimum.

\subsection{Ultra-copowers and associated maps}

Given a compact space~$Y$ and a set~$I$ we consider the \v{C}ech-Stone 
compactification $\beta(Y\times I)$, where $I$~carries the discrete topology.
There are two useful maps associated with $\beta(Y\times I)$:
the \v{C}ech-Stone extensions of the projections $\pi_Y:Y\times I\to Y$
and $\pi_I:Y\times I\to I$.
Given an ultrafilter~$u$ on~$I$ we write $Y_u=\beta\pi_I\preim(u)$
and we let $q_u=\beta\pi_Y\restr Y_u$.
In the terminology of~\cite{bankston-conf} the space~$Y_u$ is the
\emph{ultra-copower} of~$Y$ by the ultrafilter~$u$ and $q_u:Y_u\to Y$ is the
associated \emph{co-diagonal map}.
A map $f:X\to Y$ between compact spaces is \emph{co-existential} if there are
a set~$I$, an ultrafilter~$u$ on~$I$ and a map~$g:Y_u\to X$ such that
$q_u=f\circ g$.

These notions can be seen as dualizations of notions from model theory
and they offer inroads to the study of compact Hausdorff spaces by 
algebraic, and in particular lattice-theoretic, means.

\subsection{Two notions from continuum theory}

On a first-order algebraic level there is not much difference between
$Y$ and $Y_u$: they have elementarily equivalent lattice-bases for their
closed sets: the map $A\mapsto Y_u\cap\cl_\beta(A\times I)$ is an elementary embedding
of such bases.
It is therefore not unreasonable to expect that the co-diagonal 
map~$q_u$ be well-behaved.
One could expect it, for example, to be \emph{confluent}, 
which means that if $C$ is a subcontinuum of~$Y$ then every component 
of~$q_u\preim[C]$ would be mapped onto~$C$ by~$q_u$.
Certainly \emph{some} component of~$q_u\preim[C]$ is mapped onto~$C$: 
the component that contains $Y_u\cap\cl_\beta(C\times I)$ (this shows
that $q_u$~is \emph{weakly} confluent).
Intuitively there should be no difference between the components, so all
should be mapped onto~$C$.
The example below disproves this intuition.

The paper~\cite{bankston-conf} gives (references for) other reasons why
it is of interest to know whether co-diagonal and co-existential maps
are confluent.

\section{The example}

We start with the Closed Infinite Broom \cite{MR1382863}*{Example~120}
$$
B=\bigl([0,1]\times\{0\}\bigr)\cup\bigcup_{n\in\omega} H_n
$$
where $H_n=\{\orpr{t}{t/2^n}:0\le t\le1\}$ is the $n$th hair of the broom.

The range space is $B$ with the limit hair extended to have length~$2$:
$$
Y=B \cup\bigl( [1,2]\times\{0\}\bigr)
$$
We denote the extended hair $[0,2]\times\{0\}$ by $H_\omega$.

The domain of the map is $B$ with an extra hair of length~$2$ along the $y$-axis:
$$
X=B  \cup\bigl(\{0\}\times[0,2]\bigr)
$$
The map $f:X\to Y$ is the (more-or-less) obvious one:
$$
f(x,y)=\begin{cases}\orpr xy&\orpr xy\in B\\
                    \orpr{y}0& x=0
       \end{cases}
$$
Thus $f$ is the identity on~$B$ and it rotates the points on the extra hair
over~$-\frac12\pi$.

\begin{claim}
 The map $f$ is not confluent. 
\end{claim}

\begin{proof}
This is easy.
The components of the preimage of the continuum $C=[1,2]\times\{0\}$ are
the interval $\{0\}\times[1,2]$ and the singleton~$\{\orpr10\}$;
the latter does not map onto~$C$.  
\end{proof}

\begin{claim}
The map $f$ is co-existential.  
\end{claim}

\begin{proof}
We need to find an ultrafilter~$u$ and a map $g:Y_u\to X$ such 
that $f\circ g$ is the co-diagonal map $q_u:Y_u\to Y$.
In fact any free ultrafilter $u$ on~$\omega$ will do.

We define two closed subsets~$F$ and $G$ of~$Y\times\omega$ and define
$g$ on the intersections $F_u=Y_u\cap\cl_\beta F$ and $G_u=Y_u\cap\cl_\beta G$
separately.
We set
$$
F=\bigcup_{n\in\omega}\biggl(\bigcup_{k\le n}\bigl(H_k\times\{n\}\bigr)\biggr)
$$
and
$$
G=\bigcup_{n\in\omega}\biggl(
      \bigcup_{n<k\le\omega}\bigl(H_k\times\{n\}\bigr)\biggr)
$$
Note that $F\cup G=Y\times\omega$ and that $F\cap G=\{\orpr00\}\times\omega$,
so that $F_u\cup G_u=Y_u$ and $F_u\cap G_u$ consists of one point, the
(only) accumulation point of $F\cap G$ in~$Y_u$.

It is an elementary verification that $q_u[F_u]=B$ and 
$q_u[G_u]=H_\omega$.
This allows us to define $g:Y_u\to X$ by cases:
on $F_u$ we define $g$ to be just~$q_u$ and on $G_u$ we define $g=R\circ q_u$,
where $R$~rotates the plane over~$\frac12\pi$.
These definitions agree at the point in~$F_u\cap G_u$ and give continuous
maps on~$F_u$ and $G_u$ respectively.
Therefore the combined map $g:Y_u\to X$ is continuous as well.
\end{proof}

This also shows that the co-diagonal map $q_u$ is not confluent: 
no component of the preimage under~$g$ of~$\orpr10$ is mapped onto~$C$.

\begin{remark}
In \cite{bankston-interactive} 
Bankston showed that if a continuum~$K$ is such that every
co-existential map onto~$K$ is confluent then every $K$ must be 
connected im kleinen at each of its cut points.
The continuum~$Y$ above is connected im kleinen at all cut points but one:
the point~$\langle1,0\rangle$, so $Y$~does not qualify as a counterexample
to the converse.

To obtain a countereample multiply $X$ and $Y$ by the unit interval and
multiply~$f$ by the identity.
The proof that the new map is co-existential but not confluent is an easy
adaptation of the proof that $f$~has these properties.
Since $Y$ has no cut points it is connected im kleinen at all of them.
\end{remark}

\end{document}